\documentclass[11pt]{article}

\usepackage{float}
\usepackage{graphicx}%
\usepackage{multirow}%
\usepackage{amsmath,amssymb,amsfonts}%
\usepackage{amsthm}%
\usepackage{mathrsfs}%
\usepackage[title]{appendix}%
\usepackage{xcolor}%
\usepackage{textcomp}%
\usepackage{manyfoot}%
\usepackage{booktabs}%
\usepackage{algorithm}%
\usepackage{algorithmicx}%
\usepackage{algpseudocode}%
\usepackage{listings}%
\usepackage[hidelinks,bookmarksnumbered=true]{hyperref}
\usepackage{enumitem}   % to use enumerate
\usepackage{cleveref}   % to use Cref
\usepackage{bm}         % to use \bm x etc.
\usepackage{bbm}        % to use \mathbbm{1}
\usepackage{bigints}    % to use \infty
\usepackage[square,numbers]{natbib}
\bibpunct[,]{(}{)}{,}{a}{}{,}%

\usepackage[symbol]{footmisc}
\usepackage{bookmark}   % to allow paragraph to replace other section level
\usepackage{authblk}    % for author names

\usepackage{caption}                % to include some space between caption and the table
\usepackage{tikz}       % to use tikz drawing
\usetikzlibrary{calc}
\usetikzlibrary{plotmarks}
\usetikzlibrary{backgrounds}
\usetikzlibrary{arrows.meta,positioning}
\usetikzlibrary{decorations.pathreplacing}
\usetikzlibrary{patterns}
\usepackage{xcolor}

\usepackage{pgfplots}
\pgfplotsset{compat=newest}
\pgfdeclareplotmark{mystar}{
\node[star,star point ratio=2.25,minimum size=6pt,
inner sep=0pt,draw=black,solid,fill=red] {};
}
\usepackage{subfigure}
\usepackage{graphicx}

\usepackage{multirow}
\usepackage{longtable}
\usepackage{tabularx}
\usepackage{makecell}
\usepackage{array}

\usepackage{thm-restate}    % to use restatable

\def\hat{\widehat}

\def \V{\mathbf{V}}

\def\A{{\mathcal A}}

\def\C{{\mathcal C}}

\def\T{{\mathcal T}}
\def\H{{\mathcal H}}
\def\J{{\mathcal J}}

\def\T{{\mathcal T}}

\def\V{{\mathcal V}}

\def\vx{\bm{x}}
\def\vs{\bm{s}}
\def\vf{\bm{f}}

\let\cline\cmidrule

\def\rev{}

\definecolor{bluegreen}{RGB}{0,158,115}

\title{Learning \rev{Large Neighborhood Search for} Maritime Inventory Routing Optimization}
\author[1]{Rui Chen\thanks{\href{mailto:rchen@cuhk.edu.cn}{rchen@cuhk.edu.cn}}}
\author[2]{Defeng Liu\thanks{\href{mailto:defeng.liu@polymtl.ca}{defeng.liu@polymtl.ca}}}
\author[3]{Nan Jiang\thanks{\href{mailto:nanjiang@ust.hk}{nanjiang@ust.hk}}}
\author[4]{Rishabh Gupta\thanks{\href{mailto:rishabh.gupta1@exxonmobil.com}{rishabh.gupta1@exxonmobil.com}}}
\author[4]{Mustafa Kilinc\thanks{\href{mailto:mustafa.kilinc@exxonmobil.com}{mustafa.kilinc@exxonmobil.com}}}
\author[5]{Andrea Lodi\thanks{\href{mailto:al748@cornell.edu}{al748@cornell.edu}}}

\affil[1]{The Chinese University of Hong Kong, Shenzhen}
\affil[2]{Polytechnique Montr\'eal}
\affil[3]{The Hong Kong University of Science and Technology}
\affil[4]{ExxonMobil Technology and Engineering Company}
\affil[5]{Jacobs Technion-Cornell Institute, Cornell Tech and Technion - IIT}

\date{}

\begin{document}
\maketitle	

%--Paper--
\begin{abstract}

\rev{Maritime inventory routing optimization is an important yet challenging combinatorial optimization problem.} We propose a machine learning-based local search approach for finding feasible solutions of large-scale maritime inventory routing optimization problems. Given the combinatorial complexity of the problems, we integrate a graph neural network-based neighborhood selection method to enhance local search efficiency. Our approach enables a structured exploration of \rev{different} neighborhoods \rev{by imitating an optimization-based expert neighborhood selection policy}, improving solution quality while maintaining computational efficiency. Through extensive computational experiments on realistic instances, we demonstrate that our method outperforms direct mixed-integer programming \rev{as well as benchmark local search approaches} in solution time \rev{and solution quality}. 

\end{abstract}

\section{Introduction}
\label{sec_introduction}
\rev{Maritime Inventory Routing Optimization (MIRO) arises in vendor-managed supply chains where a supplier must coordinate shipping routes and inventories across multiple ports. In a MIRO problem, a fleet of vessels transports a product from production ports to consumption ports over a planning horizon, while ensuring no port’s inventory falls outside its specified bounds. This integrates two complex tasks, routing and inventory management, into a single problem. The classical Inventory Routing Optimization (IRO) problem \citep{bell1983improving} considers a single depot serving many customers with land-based vehicles over time. In contrast, MIRO typically involves multiple supply ports and multiple demand ports, with vessels that may start and end at different locations (an ``open” routing problem). Ocean-based logistics also feature much longer travel times, extended planning horizons, heterogeneous fleets, and specialized constraints (e.g., fuel or cargo-handling requirements). As noted by \cite{fagerholt2023maritime}, features such as multi-port networks, open routes, and long horizons make MIRO especially challenging compared to land-based IRO problems.} 

\rev{MIRO is of great practical importance in global trade. Maritime shipping carries roughly 9 billion tons of goods annually, accounting for over 80\% of world trade by volume. Bulk products such as crude oil, liquefied natural gas (LNG), chemicals, and cement are typically distributed over long time horizons, ranging from several months to a full year, using coordinated maritime networks. In particular, LNG supply chains are a prototypical example of long-horizon MIRO. A single liquefaction (production) terminal may need to serve multiple regasification (consumption) ports globally, with voyages lasting several weeks. Planning over a 120-360 day horizon is required to satisfy delivery contracts, manage tankfarm inventories, account for seasonal demand fluctuations, and maintain fleet utilization. These ``deep-sea” logistics problems involve large-scale networks, vessel idle times, port capacity limitations, and complex sailing schedules. In many such systems, inventory holding costs are sunk or negligible (e.g., through intra-company transfer pricing), and the objective is minimizing transportation and terminal operation costs. Given the combinatorial scale of practically relevant MIRO instances, solving them with exact mixed-integer programming (MIP) approaches becomes intractable beyond small planning horizons \citep{shaabani2023matheuristic}. As a result, practitioners increasingly rely on heuristics and matheuristics to compute high-quality solutions within reasonable time frames. This motivates the development of efficient heuristics for large-scale, long-horizon MIRO.}

\rev{In this work, we focus on the development of a graph neural network (GNN)-guided neighborhood search algorithm designed for solving long-horizon MIRO problems. While our numerical experiments use LNG supply chain planning as a representative case, given its prevalence in the literature and the real-world occurrence of long planning horizons in LNG MIRO, the proposed framework is broadly applicable to other MIRO settings with similar structural characteristics.}

\rev{The remainder of this paper is organized as follows: Section~\ref{sec_literature} reviews related literature on solution methods for MIRO and outlines the contributions of this work. Section~\ref{sec_formulation} describes the MIP formulation and modeling challenges. Our proposed learning-based approach is detailed in Sections~\ref{sec_ml_model} and~\ref{sec_algo}, followed by computational results in Section~\ref{sec_numerical}. We conclude the paper with a summary and directions for future research in Section~\ref{sec_conclusion}.}

\subsection{Literature Review}
\label{sec_literature}
A variety of solution methodologies have been proposed for the MIRO problem, ranging from exact methods to heuristics and hybrid approaches. While exact methods \rev{provide} optimality guarantees, they struggle to scale beyond small problem instances, making them impractical for real-world applications with long planning horizons. \rev{As a result, most of the literature resorts to heuristic methods rather than direct exact optimization. In this section, we provide a concise overview of existing solution approaches to show the diversity of solution strategies explored in the MIRO literature. We conclude by identifying a specific gap in neighborhood search methods, an area to which this paper aims to contribute. For comprehensive surveys of the MIRO research landscape, we refer the reader to \citet{papageorgiou2014mirplib}, \citet{fagerholt2023maritime}, and \citet{christiansen2025fifty}.}

\subsubsection{Exact Methods}
This class of  MIRO as a mixed-integer programming problem and attempted to solve it using commercial solvers, sometimes enhanced with custom cutting-plane strategies to strengthen the formulation during branch and bound \citep{engineer2012branch, andersson2017creating}. \rev{Despite these enhancements, direct solution of MIP-based formulations are generally limited to small- or medium-scale instances due to the exponential growth of the solution space \citep{shaabani2023matheuristic}.} 

To improve tractability, several papers have applied Dantzig-Wolfe and Benders decomposition methods. \citet{fagerholt2000combined} \rev{considered Dantzig–Wolfe column generation within a branch-and-bound framework for a combined ship scheduling and inventory problem. Building on this idea,} \citet{gronhaug2010branch} developed a branch-and-price method for an LNG inventory routing problem, with the master problem handling inventory and port capacity constraints while subproblems generate ship route columns. \citet{papageorgiou2014two} proposed a two-stage Benders decomposition approach that first determined \rev{high-level} inter-regional routing decisions and then optimized \rev{detailed} intra-regional loading and discharging operations.

\rev{Nonetheless, exact approaches remain difficult to scale. For instance, \citet{shao2015hybrid} showed that a Dantzig-Wolfe decomposition struggled to close the optimality gap for LNG MIRO with a 365-day horizon, underscoring the need for more scalable alternatives.}

\subsubsection{Heuristic Methods}

\rev{Given the computational challenges associated with exact methods, researchers have increasingly explored advanced heuristics and hybrid strategies for solving MIRO. One of the earliest works in this direction is that by \citet{ronen2002marine}, who introduced a ``cost-based heuristic” for a MIRO instance.}

\rev{While metaheuristics, general-purpose optimization frameworks such as genetic algorithms, particle swarm optimization, or ant colony optimization, have been explored to some extent \citep{chandra2013multi, de2017sustainable, siswanto2019maritime, al2022optimization}, purely metaheuristic approaches remain relatively few in the MIRO literature. This is largely due to the problem’s combinatorial complexity, which makes it challenging to even identify feasible solutions without domain-specific guidance or tailored feasibility checks \citep{sanghikian2025heuristic}.}

\rev{Instead, the dominant direction in the literature has focused on matheuristics that are hybrid methods, combining heuristic frameworks with mathematical programming techniques to ensure feasibility. One of the most common strategies in this category is the rolling-horizon heuristic (RHH), which decomposes long planning horizons into tractable segments. For instance, \citet{Rakke2011} applied an RHH approach to full-year LNG scheduling problems, while \cite{al2016optimal} used a similar approach for LNG vessel planning. \cite{Agra2014} explored a hybrid method for a short-sea shipping IRO (with different characteristics than LNG) that combined a RHH with local branching and a feasibility pump to obtain a good feasible solution.  More recently, \citet{li2023maritime} examined various rolling configurations and demonstrated that with proper tuning, RHH can efficiently generate high-quality solutions.}

\rev{Neighborhood-based heuristics \citep{ahuja2002survey} have also proven effective, particularly for LNG inventory routing. \citet{Stalhane2012} combined greedy insertion with neighborhood search and custom branch-and-bound procedures for a large LNG case. \citet{goel2012large} proposed a large neighborhood search (LNS) approach that utilized time-expanded arc-flow models and neighborhood operators, including a two-ship improvement heuristic (first proposed by \cite{song2013maritime}). \citet{shao2015hybrid} extended this work with a GRASP-based initialization and a broader set of local search operators, showing improved performance. }

\rev{Beyond LNG-specific applications, \citet{agra2016mip} developed neighborhood search methods for stochastic MIRO, while \citet{munguia2019tailoring} proposed a parallel LNS framework, Parallel Alternating Criteria Search, to accelerate search performance. \citet{papageorgiou2017recent} benchmarked a wide variety of matheuristic combinations, such as $k$-opt search, local branching, and rolling horizon planning, demonstrating their effectiveness on large-scale MIRO test sets. More recently, researchers have combined exact and heuristic ideas: \citet{friske2020multi} and \citet{friske2022} applied relax-and-fix and fix-and-optimize matheuristics, combining MIP-based refinement with heuristic initialization procedures to tackle challenging instances.}

\rev{Adaptive Large Neighborhood Search (ALNS) has also been applied to MIRO. ALNS uses a portfolio of neighborhood operators with adaptive weights, allowing the algorithm to ``learn” which neighborhood moves are most fruitful. \citet{hemmati2016} applied a two-phase ALNS-based matheuristic to a multi-product short-sea MIRO, initially discretizing deliveries into candidate cargoes, followed by an adaptive search over feasible cargo assignments to determine vessel schedules. \citet{heidari2023} refined this approach by integrating an ALNS with a penalized MIP that jointly optimizes vessel routing and cargo planning decisions, reporting significant performance gains.}

\subsubsection{The Neighborhood Selection Bottleneck}

\rev{The literature reviewed in the previous subsection makes clear that heuristic approaches offer the best scalability for large MIRO instances. In this paper, we focus specifically on neighborhood-based heuristics, which have shown particular promise for long-horizon MIRO instances \citep{goel2012large, shao2015hybrid, papageorgiou2017recent}}.

\rev{The effectiveness of LNS depends heavily on the strategy used to select the neighborhood to ``destroy and repair” in each iteration. This selection step involves a trade-off between speed and informativeness. Random selection is computationally trivial but blind to problem structure, often leading to wasted iterations and slow convergence. At the other extreme, oracle-based strategies evaluate many neighborhoods (e.g., via LP relaxations) to identify the most promising one, but their computational cost is often prohibitive in practice.}

\rev{Adaptive methods like ALNS aim to balance this trade-off by learning which operators have performed well historically. However, these strategies only adapt based on past performance and lack predictive insight into the current problem state. ALNS cannot leverage the structural information embedded in vessel schedules, port congestion, or inventory constraints to guide selection effectively.}

\rev{This highlights an important gap in the existing literature: current neighborhood selection strategies are either too simplistic to be effective or too computationally expensive for use in real-time applications. There is a clear need for methods that can make informed, structure-aware neighborhood selections while maintaining computational efficiency. Machine learning offers a promising avenue in this regard, as it is well suited to learning mappings from high-dimensional problem states to effective neighborhood selection strategies.}

In this work, we propose using machine learning (ML) to guide the selection of neighborhoods in an LNS scheme designed to solve large-scale instances of the LNG-IRO problem. By learning patterns from high-quality solutions, our approach aims to identify promising neighborhoods more effectively, reducing unnecessary evaluations and improving the overall efficiency of the search process. This represents a novel direction in the application of machine learning to large-scale optimization problems, addressing a key limitation of existing LNS methods.

\subsubsection{ML-Augmented Local Search}
\label{sec_ml}

The last decade has witnessed an impressive amount of work devoted to exploring the synergies between machine learning and discrete optimization. Those synergies go both ways, with combinatorial optimization (CO) applied to many significant learning problems, for example in the area of interpretable ML \citep{rudin2018learning, rudin2022interpretable}, and with ML applied in a variety of ways to help solving of combinatorial optimization problems \citep{bengio2021machine, KFVHW21, scavuzzo2024machine}.

For the latter case, the first attempted paradigm has been ``end-to-end" learning, that is constructing a heuristic solution directly by training an ML model \citep{khalil2017learning_graph, nazari2018reinforcement, kool2018attention, zhang2020learning, liu2021learning}. Those methods are typically designed for specific CO problems that are not (tightly) constrained, and a heuristic solution can be easily constructed. Conversely, since a wide range of constrained CO problems can be formulated into a mixed-integer programming (MIP) model, there has also been increasing interest in learning decision rules to improve MIP algorithms \citep{gasse2019exact, sonnerat2021learning, nair2020solving, bonami2022classifier}. The role of graph neural networks (GNNs) \citep{gori2005new,scarselli2008graph, hamilton2017representation, cappart2023combinatorial} in developing ML methods for CO has been impressive, \rev{starting with the GNN-based branching scheme in \cite{gasse2019exact}, the first to report computational results competitive with the state of the art in integer programming}.

Especially relevant in the context of the current paper are the efforts to improve local search, especially LNS algorithms for combinatorial optimization, see, e.g., \citet{liu2023machine, Liu_Fischetti_Lodi_2022}. There, the scope is to learn structural properties of the neighborhoods that are iteratively explored by local search and metaheuristic methods to compute high-quality solutions in a reasonable amount of computing time, that is, when solving the problem of interest to proven optimality is out of reach. Such learning might be crucial for multiple reasons, the most common one being that the local search scheme can pick a ``good" neighborhood without resorting to computationally expensive analyses, for example by learning to imitate an expensive expert. 

This is precisely what we are proposing to do here: It is computationally  expensive to find even a single feasible solution for the MIRO formulation discussed in the next section, so we construct an associated feasibility problem and we solve it by local search, exploring large neighborhoods learned by sophisticated ML methods. The approach is detailed in Sections \ref{sec_ml_model} and \ref{sec_algo}.

\rev{The primary contributions of this work are as follows:}
\begin{enumerate} 
\item \rev{A Novel ML-Guided LNS Framework for MIRO: We propose and implement a machine learning framework to guide neighborhood selection specifically for large-scale MIRO. We formalize the neighborhood selection task as a policy learning problem and demonstrate how to train a model to imitate a strong, optimization-based expert policy (the Root policy), effectively distilling its expensive computational wisdom into a fast, deployable neural network.}

\item \rev{A GNN-Based Neighborhood Scoring Model: We design a specific GNN architecture that effectively encodes the complex, multipartite structure of the MIRO, including vessels, ports, tankfarms, and contracts, into a unified graph representation. This model learns to process the global state of the current solution and predicts the potential of each two-vessel neighborhood, enabling an intelligent and targeted search that focuses computational effort where it is most likely to yield improvements.}

\item \rev{An Effective Imitation Learning Strategy: We employ the Dataset Aggregation (DAGGER) algorithm for training our GNN policy. This iterative imitation learning approach is demonstrably more effective than standard offline supervised learning for sequential decision-making problems. By collecting training data from trajectories generated by the learned policy itself, DAGGER mitigates distributional shift and improves the policy's performance and stability in the online search setting.}

\item \rev{Extensive Empirical Validation: We conduct extensive computational experiments on a suite of realistic, large-scale MIRO instances designed to mirror real-world industrial challenges. We demonstrate that our learning-based LNS approach significantly outperforms both direct MIP solves with a state-of-the-art commercial solver and benchmark LNS approaches in finding feasible solutions. Specifically, our method reduces total computing time by up to 40.48\% on instances with purely combinatorial constraints and 13.94\% on hybrid instances that also include inventory constraints.}
\end{enumerate}

\section{Problem Formulation}
\label{sec_formulation}
\rev{In this section, we present the core formulation of the MIRO problem considered in this work. We start by providing an overview of MIRO in the LNG industry to motivate and contextualize our modeling choices. Many of the complexities observed in LNG logistics are also present in other maritime inventory routing settings. Thus, while LNG serves as the central motivating application, the formulation and solution methodology are broadly applicable. The core model we describe builds upon and extends existing MIRO formulations in the literature (e.g., \citealp{song2013maritime, papageorgiou2014mirplib}).}

\subsection{\rev{LNG Industry Context}}
\rev{The global LNG industry has become increasingly complex as major players shift from long-term bilateral supply agreements toward more dynamic, portfolio-based trading models. In this evolving landscape, companies engage in both upstream production and downstream sales, while actively participating in the spot market to optimize margins. This growing commercial flexibility introduces significant operational complexity, requiring careful coordination of supply, inventory, and transportation assets.}

\rev{A typical LNG portfolio includes a mix of internal production sources, coupled with limited storage capacity in the form of tankfarms and external supply via long-term \textit{purchase contracts}. These contracts commonly include \textit{ratability constraints}, which restrict the number and timing of LNG cargoes that can be purchased within specified time windows. On the demand side, \textit{sales contracts} impose similar constraints on delivery timing and quantities. Additionally, contracts often specify \textit{incoterms}, which determine which party is responsible for transportation. In cases where the seller is responsible for shipping, it becomes essential to coordinate vessel schedules to fulfill contractual obligations efficiently.}

\rev{Shipping logistics add another layer of complexity. LNG cargoes are transported using a limited fleet of heterogeneous vessels, comprising both company-owned ships and higher-cost spot-chartered vessels. Each vessel operation must respect several practical constraints, including berth availability (only one vessel can occupy a berth at a time), port turnaround times, compatibility with LNG types, and maintenance or bunkering requirements. Importantly, vessel utilization decisions must balance cost-efficiency and contract fulfillment, as chartering additional vessels significantly increases operational costs.}

\rev{The LNG MIRO is inherently both a planning problem and a scheduling problem. In scenarios involving multiple supply sources and customers, as considered in this paper, the core difficulty lies in determining a set of interdependent decisions that collectively satisfy all logistical, operational, and contractual constraints. These decisions can be broadly categorized as follows:}

\begin{enumerate}
\item \rev{\textbf{Cargo planning:} Select which supply sources (tankfarms or purchase contracts) fulfill which demand contracts, and specify the timing of each cargo’s loading and delivery. This decision must respect time window constraints, ratability requirements, and tankfarm inventory limits. Even without optimizing economic objectives, identifying a feasible assignment of cargoes that satisfies all constraints is a highly combinatorial task.}
\item \rev{\textbf{Vessel scheduling:} Determine the routing and timing of each vessel’s operations, ensuring compatibility with port capacities, sailing durations, loading/unloading constraints, and contract fulfillment. The complexity of coordinating these activities across a heterogeneous fleet and over a long planning horizon makes the generation of feasible schedules computationally challenging.}
\end{enumerate}

\subsection{\texorpdfstring{\rev{Modeling Assumptions}}{Modeling Assumptions}}
\rev{We make the following standard modeling assumptions:}

\begin{enumerate}
\item \rev{Each tankfarm and contract is associated with a unique port.}
\item \rev{Travel times between ports are deterministic.}
\item \rev{Production rates at each tankfarms are pre-determined.}
\item \rev{Inventory at each tankfarm must remain within specified bounds throughout the planning horizon.}
\item \rev{Each vessel is always loaded to its full capacity at the supply source; partial loading is not considered.}
\item \rev{A vessel performs exactly one loading and one delivery operation per voyage; there is no transshipment or multi-drop routing.}
\item \rev{Ship-to-ship transfers, weather disruptions, and canal congestion are not explicitly modeled.}
\end{enumerate}

\subsection{\rev{Modeling Approach}}
\begin{figure}[tp!]
\includegraphics[width=\textwidth]{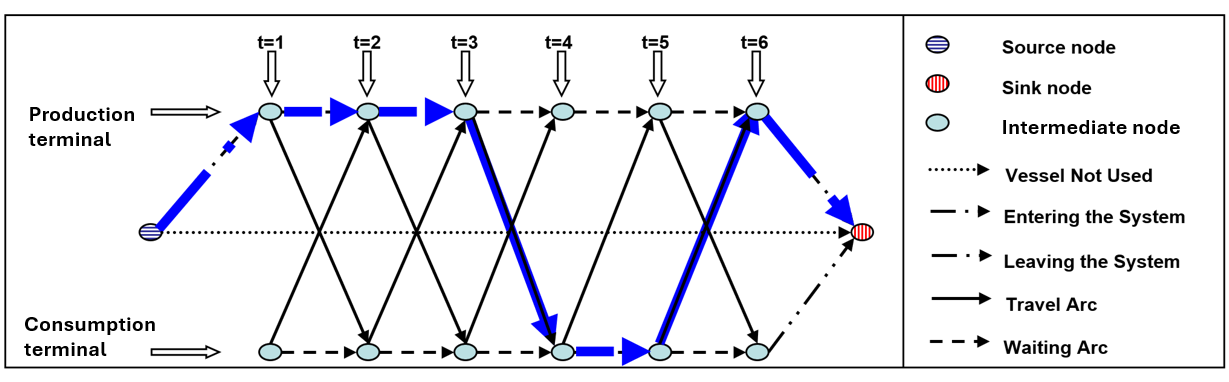}
\caption{A TSN with a single production port and a single consumption port}
\label{figure:TSN}
\end{figure}

\rev{To model the complex spatial and temporal dynamics of this real-world system, we adopt a Time-Space Network (TSN) formulation, similar to that introduced by \cite{song2013maritime}. The planning horizon is discretized into uniform time intervals (e.g., daily periods), and the movement and operations of each vessel are captured in a directed graph. As shown in Figure \ref{figure:TSN}, for each vessel $v \in \mathcal{V}$, we define a graph $\mathcal{G}_v = (\mathcal{N}_v, \mathcal{A}_v)$, where each node $n \in \mathcal{N}_v$ represents a specific port at a specific time, and each arc $a \in \mathcal{A}_v$ denotes a feasible action, such as sailing between ports or waiting at a location. The vessel schedule corresponds to a path through its TSN, beginning at a source node (representing its initial position and start time) and ending at a sink node (representing the end of the planning horizon). Nodes on the row of ``production terminal'' represent where and when the vessel gets loaded, while on the row of ``consumption terminal'' represent where and when the cargoes get delivered.}

\rev{The problem input data can be broadly categorized into individual attributes and pairwise attributes. Individual attributes include characteristics of}
\begin{itemize}
\item \rev{Vessels (e.g., type, cargo capacity, initial location),}
\item \rev{Ports (e.g., berth availability, loading/discharging durations),}
\item \rev{Tankfarms (e.g., production rates, inventory bounds), and}
\item \rev{Contracts (e.g., time windows, origin/destination ports, incoterms).}
\end{itemize}

\rev{Pairwise attributes define relationships between components, such as}
\begin{itemize}
\item \rev{Sailing times and distances between port pairs,}
\item \rev{Compatibility constraints between vessels and contracts, and}
\item \rev{Feasible loading/discharging operations at specific port-vessel combinations.}
\end{itemize}

\rev{This TSN-based modeling approach enables the explicit representation of temporal and spatial dynamics inherent in LNG logistics. It allows integration of contract economics, vessel operations, and inventory management, providing a framework for solving the LNG MIRO problem at both strategic and operational levels.}  

\subsection{\rev{Core Mathematical Formulation}}
\rev{In this work, \textbf{our objective is to compute a feasible vessel schedule} that satisfies key operational and contractual constraints. We formulate this as a MIP over a TSN representation. Here, we present a core formulation that captures the essential structure of the LNG MIRO problem. Additional real-world complexities, considered in our actual implementation and modeled through supplementary variables and constraints, are omitted here for brevity. A brief discussion of these additional considerations is provided following the core formulation.}

\rev{
\paragraph{Sets and Indices:}
\begin{enumerate}
\item $\mathcal{H}=\{1,\ldots,|\mathcal{H}|\}$ - set of discrete time-periods in the planning horizon (index $h$).
\item $\mathcal{V}$ - set of vessels (index $v$). 
\item $\mathcal{T}$ - set of tankfarms from where product can be loaded (index $t$).
\item $\mathcal{C}$ - set of purchase and sales contracts (index $c$). All purchase contracts are represented by the set $\mathcal{C}_l$ and sales contracts by $\mathcal{C}_u$ such that $\mathcal{C} = \mathcal{C}_l \cup \mathcal{C}_u$. 
\item $\mathcal{J}$ - set of all port locations (index $j$). All loading sources (tankfarms and purchase contracts ports) are represented by the set $\mathcal{J}_l$ and delivery ports by $\mathcal{J}_u$. 
\item $\mathcal{N}_v$ - set of TSN nodes (index $n$) for each vessel $v \in \mathcal{V}$, define $\mathcal{N}_v = \{(j,h) : j\in \mathcal{J},\ h\in \mathcal{H}\}$ as the nodes representing vessel $v$ being at port $j$ at time $h$.
\item $\mathcal{A}_v$ - set of directed arcs (index $a$) in the time-space network for vessel $v$. These arcs consist of either sailing arcs (denoted by set $\mathcal{A}_v^\text{sailing}$) of the form $((j,h), (j', h+T_{jj'}))$ representing vessel $v$ traveling from port $j$ to port $j'$ (with $T_{jj'}$ travel time), or waiting arcs ($\mathcal{A}_v^\text{waiting}$) of the form $((j,h), (j, h+1))$ representing vessel $v$ staying idle at port $j$ for one period. We denote by $\delta^-_{(j,h)}$ the set of arcs (both sailing and waiting) entering node $(j,h)$ and by $\delta^+_{(j,h)}$ the set of arcs leaving node $(j,h)$, for a given vessel’s TSN.
\item $\mathcal{R}_c$ - set containing contractual ratability (index $r$) requirements for contracts $c \in \mathcal{C}$. This set contains tuples $(h_1, h_2, LB, UB)$ specifying that at least $LB \in \mathbb{Z}_+$ and at most $UB \in \mathbb{Z}_{+}$ cargoes must be delivered (or loaded) between times $h_1 \in \mathcal{H}$ and $h_2 \in \mathcal{H}$ (inclusive). 
\end{enumerate}}

\rev{\paragraph{Parameters:}
\begin{enumerate}
\item $S_v \in \mathcal{J}$ - initial location of vessel $v \in \mathcal{V}.$
\item $E_v \in \mathcal{J}$ - required end location of vessel $v \in \mathcal{V}.$
\item $I^{\mathrm{initial}}_{t} \in \mathbb{R}_+$ - initial inventory level at tankfarm $t\in\T$ at the start of the time horizon.
\item $\bar{I}_t \in \mathbb{R}_+$, $\underline{I}_t \in \mathbb{R}_+$ - inventory upper and lower bounds at tankfarms $t \in \mathcal{T}$ (capacity and safety stock limits, respectively).
\item $P_{th} \in \mathbb{R}_+$ - production amount at tankfarm $t\in\T$ in period $h\in\H$.
\item $C_v \in \mathbb{R}_+$ - capacity of vessel $v\in\V$.
\item $B_j \in \mathbb{Z}_+$ - berth capacity at port $j\in\J_l$ limiting the number of vessels that can simultaneously (un)load at a port $j$.
\item $T_{jj'} \in \mathbb{Z}_+$ - travel time (in periods) for a vessel to sail from port $j$ to port $j'$.
\item $J_t\in\J$ - the unique loading port $j\in\J$ associated with tankfarm $t\in\T$.
\item $J_c\in\J$ - the unique (loading/discharging) port $j\in\J$ associated with (purchase/sales) contract $c\in\C$.
\end{enumerate}}

\rev{
\paragraph{Decision Variables:}
\begin{enumerate}
\item $x_{va} \in \{0,1\}$ - binary variable for vessel $v\in \V$ and arc $a\in \A_v$. Specifically, $x_{va} = 1$ if and only if vessel $v$ traverses arc $a$ in the schedule (i.e., travels or waits from one port-time node to another along that arc). These variables define each vessel’s route and schedule over the time-space network.
\item $z^l_{vjh} \in \{0,1\}$ - binary loading indicator for vessel $v\in\V$ and port $j\in \mathcal{J}_l$ at time $h\in \H$. Specifically, $z^l_{vjh}=1$ if and only if vessel $v$ loads the product at port $j$ in time-period $h$.
\item $z^u_{vjh} \in \{0,1\}$ - binary discharging indicator for vessel $v$ at port $j\in \J_u$ at time $h\in \H$. Specifically, $z^u_{vjh}=1$ if vessel $v$ discharges the product at port $j$ in time-period $h$.
\item $i_{th} \ge 0$ - inventory level at tankfarm $t \in \mathcal{T}$ at the end of period $h$.
\end{enumerate}
Using the above notation, we formulate the LNG MIRO as a mixed-integer feasibility problem
%under vessel routing constraints, inventory flow conservation, and operational constraints 
as detailed below.}

\rev{
\paragraph{Constraints:}
\begin{enumerate}
\item {\textbf{Vessel Routing Constraints}}\\
Each vessel must depart from its specified initial location at the start of the planning horizon:
\begin{align*}
\sum_{a \in \delta^+_{(S_v,1)}} x_{va} = 1 \quad \forall v \in \mathcal{V}.
\end{align*}
Each vessel must arrive at its designated end location at the end of the planning horizon:
\begin{align*}
\sum_{a \in \delta^-_{(E_v,|\mathcal{H}|)}} x_{va} = 1 \quad \forall v \in \mathcal{V}.
\end{align*}
Flow conservation must be maintained at every node in the network:
\begin{align*}
\sum_{a \in \delta^+_{(j,h)}} x_{va} = \sum_{a \in \delta^-_{(j,h)}} x_{va} 
\quad \forall v \in \mathcal{V},\, n=(j,h)\in \mathcal{N}_v\setminus\{(S_v,1),(E_v,|\mathcal{H}|)\}.
\end{align*}
\item {\textbf{Loading and Delivery Constraints}}\\
A vessel is considered to have loaded cargo at a supply port if it departs that port on a sailing arc:
\begin{align*}
z^l_{vjh} = \sum_{a \in \mathcal{A}^{\mathrm{sailing}}_v \cap \delta^+_{(j,h)} } x_{va}
\quad \forall v \in \mathcal{V},\, j \in \mathcal{J}_l,\, h \in \mathcal{H}.
\end{align*}
A vessel is considered to have delivered cargo at a delivery port if it arrives there on a sailing arc:
\begin{align*}
z^u_{vjh} = \sum_{a \in \mathcal{A}^{\mathrm{sailing}}_v \cap \delta^-_{(j,h)} } x_{va}
\quad \forall v \in \mathcal{V},\, j \in \mathcal{J}_u,\, h \in \mathcal{H}.
\end{align*}
\item {\textbf{Inventory Constraints}}\\
Inventory levels at each tankfarm start from a specified initial value:
\begin{equation}\label{constr:initial_inventory}
i_{t1} = I^{\mathrm{initial}}_t \quad \forall t \in \mathcal{T}.
\end{equation}
Inventory evolves over time based on production and cargo loadings:
\begin{equation}\label{constr:inventory_change}
i_{t(h+1)} = i_{th} + P_{th} - \sum_{v \in \mathcal{V}} C_v z^l_{vJ_th}
\quad \forall t \in \mathcal{T},\, h \in \mathcal{H}.
\end{equation}
Inventory levels must remain within storage capacity and safety stock limits:
\begin{align*}
\underline{I}_t \leq i_{th} \leq \overline{I}_t \quad \forall t \in \mathcal{T},\, h \in \mathcal{H}.
\end{align*}
\item {\textbf{Combinatorial (Port Capacity and Ratability) Constraints}}\\
The number of vessels loading at any port in a given time period cannot exceed the berth capacity:
\begin{align*}
\sum_{v \in \mathcal{V}} z^l_{vjh} \leq B_j \quad \forall j \in \mathcal{J}_l,\, h \in \mathcal{H}.
\end{align*}
The number of vessels discharging at a delivery port must also respect berth limits:
\begin{align*}
\sum_{v \in \mathcal{V}} z^u_{vjh} \leq B_j \quad \forall j \in \mathcal{J}_u,\, h \in \mathcal{H}.
\end{align*}
Contract ratability constraints require a minimum and maximum number of cargoes to be purchased (or delivered) within specified time windows for each contract:
\begin{align*}
LB \leq \sum_{h=h_1}^{h_2} \sum_{v\in\V}z^l_{vJ_ch} \leq UB 
\quad \forall c \in \mathcal{C}_l,\, (h_1,h_2,LB,UB) \in \mathcal{R}_c,\\
LB \leq \sum_{h=h_1}^{h_2} \sum_{v\in\V} z^u_{vJ_ch} \leq UB 
\quad \forall c \in \mathcal{C}_u,\, (h_1,h_2,LB,UB) \in \mathcal{R}_c.
\end{align*}
\end{enumerate}
}

\rev{\paragraph{Remark:} While the core formulation above provides a clean abstraction of the problem, the implementation that we used to conduct the computational experiments reported in this paper requires additional constraints and variables to capture operational complexities. In our implementation, we incorporate the following extensions:
\begin{enumerate}
\item Certain vessels require bunkering (refueling) between voyages, typically after each round trip. We model these by increasing the length of the leg where the vessel refuels.
\item Some vessels are subject to maintenance windows that must be honored by routing them to designated ports within specified time windows.
\item A limited pool of spot-chartered vessels is available to handle temporary shortage of vessels. Unlike owned vessels, these chartered vessels are exempt from network flow conservation constraints; their availability can be initiated at any port in the planning horizon, and their final position is unconstrained, allowing termination at any port.
\item \rev{Contracts may specify incoterms such as Free on Board (FOB) or Delivered Ex-Ship (DES), which determine the party responsible for transportation. Under FOB terms, the buyer provides the vessel and arranges transportation, whereas under DES terms, this responsibility lies with the seller. In such cases, the optimization may only need to match supply with demand, as vessel routing decisions are handled externally. However, these externally arranged vessel operations still occupy berth capacity at the loading or unloading ports and must be accounted for in the overall scheduling to avoid conflicts and ensure operational feasibility.}
\item Each vessel becomes available for routing only after a specified release date, rather than being available from the start of the planning horizon.
\item Loading/discharging operations may span multiple time periods. While the core model assumes instantaneous transfers, our full implementation includes multi-day loading/discharging windows and berth occupation durations.
\item There can be multiple types of LNG with different properties (e.g., energy density), and certain contracts are only compatible with LNG with certain grades.
\end{enumerate}
}

Since our real-world MIRO instances are too large for commercial solvers to compute a feasible solution within a reasonable amount of time, we consider \rev{a reformulated feasibility problem. Specifically, we reformulate the feasibility problem of our MIRO problem as an infeasibility minimization problem by first adding slack variables and then minimizing the slack in the objective, which results in an optimization problem of the following form:
\begin{equation}\label{fea_MIP}
\begin{aligned}
\min_{\vx,\vs}~&\sum_{i\in I}s_i\\
\text{s.t. }&\vx_v\in \Omega_v, &&v\in \V\cup\{0\},\\
&\sum_{v\in \V\cup\{0\}}\vf_{iv}^\top \vx_v\leq r_i+\rho_i s_i, &&i\in I,\\
&s_i\geq 0, &&i\in I.
\end{aligned}
\end{equation}
Here, $\vx_v$ denotes decision variables associated with vessel $v\in \V$ (e.g., having $v$ in its index in the core formulation). Some operations do not involve vessel scheduling, e.g., DES purchase and FOB sales operations (which do not exist in the core model though), and $\vx_0$ denotes variables associated with those non-vessel operations. We consider all non-vessel operations, as a special type of vessel. Set $\Omega_v$ describes the set of vessel schedules that meet vessel constraints associated with vessel $v$, e.g., vessel routing constraints. The objective of Problem \eqref{fea_MIP} is to schedule different vessels to meet all constraints (or find a solution with the ``minimum violation"). Set $I$ denotes the set of joint constraints that involve multiple vessels, namely, inventory constraints, port capacity and ratability constraints. For each such constraint $i\in I$ and vessel $v\in \V\cup\{0\}$, $\vf_{iv}$ denotes how variables associated with vessel $v$ contribute the $i$-th joint constraint, and $s_i$ denotes the slack variable (scaled by weight $\rho_i>0$) associated with joint constraint $i$. Scaling parameters $\rho_i$ are set to $1$ for port capacity and ratability constraints, and set to the safety stock/storage capacity limit for inventory (lower/upper bounding) constraints. Also, note that there are no additional auxiliary variables (other than the slack variables) in Problem \eqref{fea_MIP} as one can project out those auxiliary variables ($i_{th}$ variables) in the core formulation using their defining constraints (i.e., equations \eqref{constr:initial_inventory} and \eqref{constr:inventory_change}). Another nice property about Problem \eqref{fea_MIP} is that it is very simple to construct a feasible solution by picking a particular solution from $\Omega_v$ for each vessel $v\in \V\cup\{0\}$ and setting the values of the slack variables accordingly.
}

\rev{In this paper, we consider a large neighborhood search approach to tackle problem \eqref{fea_MIP}.}
In particular, we focus on the ``two-ship improvement heuristic'' as proposed in \cite{goel2012large}, which has the most promising performance in our preliminary experiments. We will specify our large neighborhood search strategies later.

\section{ML Models for Finding Good Neighborhoods}
\label{sec_ml_model}
\rev{Local search is a versatile tool for many combinatorial optimization (CO) problems \citep{aarts2018local}. The idea is to improve the current solution and restrict the search for new solutions in a lower dimensional space so that restricted subproblems can be solved efficiently. In the literature, its application to decomposable problems is especially effective and widespread. Decomposable optimization problems can be broken down into smaller localized subproblems. The nature of this structure is highly advantageous for efficient neighborhood exploration. For example, in the context of the vehicle routing problem, which is often decomposed by vehicle, local search operators that focus on individual routes have been shown to be highly effective \citep{fisher1981generalized, toth2002vehicle}. Moreover, the structure of a decomposable problem often lends itself to the design of highly effective, problem-specific neighborhood operators \citep{bandelt2004local, liu2023machine}. In the context of LNG problems, local search heuristics also have been used for solving LNG inventory routing \citep{goel2012large,shao2015hybrid}. }

Specifically, local search moves from one solution to another by applying local changes to a set of free variables. For an optimization problem of the form $\max_{\vx}\{f(\vx):\vx\in S\subseteq \mathbb{R}^n\}$, at a particular iteration with the incumbent solution $\vx=\hat{\vx}$, we choose a subset $F\subseteq[n]$ of the variable indices and solve the problem $\max_{\vx}\{f(\vx):\vx\in S;~x_i=\hat{x}_i,~i\in [n]\setminus F\}$ to obtain a new solution. This local search step is repeated until a satisfactory solution is found. A key decision in this local search process is the choice of the index set $F$ in each iteration.

\subsection{Vessel Neighborhoods}
For combinatorial optimization problems with specific structures like LNG inventory routing, the index set \rev{$F$} in each local search step is often selected based on neighborhoods defined by specific components of the problem. In the context of LNG inventory routing, such neighborhoods can be defined by vessels, time windows, contracts, or locations. For simplicity, we considered only neighborhoods defined by vessels. However, the techniques we developed are not specific to vessel neighborhoods and can be used for other more general neighborhoods as well.

Building on the two-ship improvement heuristic proposed in \cite{goel2012large}, we specifically tested local search approaches using two-vessel neighborhoods. \rev{The basic idea is to start from some initial solution of the problem, and then iteratively pick two particular ones, fix all vessel variables (slack variables are not fixed) associated with vessels other than those two picked ones, and solve the MIP subproblem with the remaining variables to update the solution. The MIP subproblem generated from each local search iteration is solved using a commercial solver with a given computational budget.} In this approach, the free variables \rev{$\vx_F$} are defined by \rev{$(\vx_v,\vx_{v'},\vs)$} associated with a pair $\{v,v'\}\subseteq \V\cup\{0\}$ with $v\neq v'$ \rev{in Problem \eqref{fea_MIP}}.
% one-vessel neighborhoods and two-vessel neighborhoods. 
% In one-vessel neighborhood search, the free variables $x_I$ are defined by variables $x_v$ associated with one $v\in V\cup\{0\}$.
% In two-vessel neighborhood search, the free variables $x_I$ are defined by $(x_v,x_{v'})$ associated with a pair $\{v,v'\}\subseteq V\cup\{0\}$ with $v\neq v'$.
% Preliminary experiments suggest that employing only one-vessel neighborhoods in the local search to solve the feasibility Problem \eqref{fea_MIP} can lead to local optima with significant feasibility gaps. However, if we apply two-vessel neighborhoods in the local search, even with randomly selected $\{v,v'\}$ from $V\cup\{0\}$, the feasibility gap would gradually shrink to $0$ and in the end a feasible solution can often be obtained.

\subsection{Neighborhood Selection Rules and the Chosen Expert}
\label{sec_neighborhood}
Since the number of \rev{choices for vessel pairs} that one can potentially use in a local search is large, finding the optimal pair of vessels to use is non-trivial (even the definition of ``\rev{optimal}" is not trivial in this case).
We tested a variety of different neighborhood selection rules, namely
\begin{enumerate}
\item \emph{Random}: uniformly sample a pair $\{v,v'\}$ from $\V\cup\{0\}$;
\item \emph{LP}: solve the linear programming relaxations \rev{(of the local search MIP subproblems)} associated with all possible pairs $\{v,v'\}$ \rev{and rank different neighborhoods according to their LP bounds};%select the one that leads to \rev{the} best dual bound;
\item \emph{Root}: solve the \rev{local search MIP subproblems} associated with all possible pairs $\{v,v'\}$ without branching \rev{(but with solver presolve and cutting planes)}, i.e., at the root node \rev{(by setting ``\texttt{NodeLimit}'' to 1 in Gurobi)}, and \rev{rank different neighborhoods according to their root bounds \rev{(i.e., the dual bounds obtained at the root node)}};% select the one that leads to the best dual bound;
\item \emph{Exact}: solve the \rev{MIP subproblems} associated with all possible pairs $\{v,v'\}$ until hitting the time limit, and \rev{rank different neighborhoods according to their primal bounds}.%select the one that leads to the best primal bound.
\end{enumerate}
Based on several preliminary experiments comparing the performance of these neighborhood selection rules in the context of LNG inventory routing, as one might expect, \emph{Random}$>$\emph{LP}$>$\emph{Root}$>$\emph{Exact} in terms of the running time. On the other hand, \emph{Random}$<$\emph{LP}$<$\emph{Root}$<$\emph{Exact} in terms of the quality of the selected neighborhoods, and the difference in quality can be significant. For example, for one particular test instance, \emph{Random} takes $188$ iterations to find a completely feasible schedule while \emph{Exact} takes only $25$ iterations. Therefore, it is evident that local search can be significantly accelerated by using strong neighborhood selection rules. However, directly applying \emph{LP}, \emph{Root}, or \emph{Exact} in our solution framework would be computationally infeasible as they are very computationally prohibitive (often orders of magnitude slower than applying \emph{Random}). We then explore the use of machine learning tools to mimic the behavior of these strong neighborhood selection rules.  

\rev{To train ML models in the imitation learning framework, we need to select an expert policy to provide the labeled demonstrations in neighborhood search. In principle, it is desirable to have a fast expert to allow for efficient exploration and data generation, but still good enough to be effective in exploiting the neighborhoods for finding improved solutions. To balance the trade-off between running-time and solution quality, we select \emph{Root} as the expert. The \emph{Root} policy has the capacity to select solution neighborhoods with good quality, and the amount of computation required for data collection is still controllable as the method only solves the root node of the MIP subproblem.}

\subsection{Graph Neural Networks for Neighborhood Search}
\label{sec_neighborhood_gnn}

Recently, different works (see, e.g., the review paper by \citet{cappart2023combinatorial}) have explored the use of GNNs in the context of combinatorial optimization. As a learning framework, GNNs have many favorable properties, such as permutation invariance, sparsity, and scalability. In the construction of the input graph, key elements in CO are encoded by nodes with (multi-dimensional) attributes, and the relationships between different elements are encoded by edges with (multi-dimensional) attributes. GNNs then iteratively update the nodes by aggregating information from their neighbors with the goal of minimizing a particular loss function.

In the context of LNG inventory routing, a (potentially infeasible) schedule can be represented as a graph with embedded nodes and edges and further encoded by a GNN model. The goal is to train a GNN model that takes an LNG schedule graph as input and predicts high scores for promising pairs of vessel nodes selected by the expert policy.
\begin{table}[tb!]
\centering
\caption{Description of the features in the LNG graph $\mathbf{g}(\mathbf{N}, \mathbf{E})$.}
% \vspace{10pt}
\label{table:lng-graph-features}
\begin{tabular}{c l l}
\multicolumn{1}{c}{} & \multicolumn{1}{c}{} & \multicolumn{1}{c}{\rev{Description of Features}} \\
\toprule
\multirow{1}{*}{Node} 
& Vessel Node & vessel type, capacity, maintenance profile, costs \\
\cmidrule{2-3}
& Tankfarm Node & safety tank top, safety tank bottom, production profile \\
\cmidrule{2-3}
& Port Node & discharging duration of the port \\
\cmidrule{2-3}
& Contract Node & contract profile with start and end dates \\
\cmidrule{2-3}
& Ratability Node & ratability constraint for each contract \\

\midrule
\multirow{1}{*}{Edge}
& Vessel-Tankfarm &  if the vessel is compatible
with the tankfarm \\
\cmidrule{2-3}
& Vessel-Port & if the vessel is compatible
with the port \\
% \midrule
% \multirow{1}{*}{$\mathbf{Index}$}
\cmidrule{2-3}
& Vessel-Contract & if the vessel is compatible with the contract \\
\cmidrule{2-3}
& Vessel-Ratability & the vessel’s contribution to the ratability \\
\bottomrule
\end{tabular}
\end{table}

\subsubsection{Feature Engineering}

We represent the information available to select neighborhoods as a multipartite graph with attributes. These attributes are very high-dimensional in comparison with the number of data points we have, as many features are indexed by combinations of vessels, ports and dates. Therefore, some feature engineering is critical here to avoid overfitting the model. \rev{Specifically, we select five types of nodes: vessel node, tankfarm node, port node, contract node and ratability node. To connect those nodes into a multipartite graph, we introduce 4 types of edges: vessel-tankfarm, vessel-port, vessel-contract and vessel-ratability. The LNG graph is denoted as $\mathbf{g}(\mathbf{N}, \mathbf{E})$, where $\mathbf{N}$ denotes the set of nodes and $\mathbf{E}$ denotes the set of edges in the graph. The features of nodes and edges are listed in Table \ref{table:lng-graph-features}.}

\subsubsection{GNN Architecture}

Our GNN design \rev{consists of three modules: the input module, the convolution module, and the output module. The input module embeds the input graph, where features of the state nodes and edges are encoded by a  Multi-Layer Perceptron (MLP) neural networks, into the latent space. The convolution module propagates the encoded features in the graph with graph convolution layers. }

\rev{At each convolution layer, our GNN architecture applies the message-passing operator, given an LNG graph $(\mathbf{N}, \mathbf{E})$ with its node set $\mathbf{N}$ and edge set $\mathbf{E}$, the message-passing computation for any node $i \in \mathbf{N}$ can be defined as
\begin{align}
    \mathbf{v}_i^{(h)} = f_{\theta}^{(h)} \left( \mathbf{v}_i^{(h-1)}, \sum_{j \in \mathcal{N}(i)} \, g_{\phi}^{(h)}\left(\mathbf{v}_i^{(h-1)}, \mathbf{v}_j^{(h-1)},\mathbf{e}_{j,i}\right) \right),
\end{align}
where $\mathbf{v}^{(h-1)}_i \in \mathbb{R}^d$ denotes the feature vector of node $i$ from layer $(h-1)$, $\mathbf{e}_{j,i} \in \mathbb{R}^m$ denotes the feature vector of edge $(j, i)$ from node $j$ to node $i$ of layer $(h-1)$, $\mathbf{f}_{\theta}^{(h)}$ and $ \mathbf{g}_{\phi}^{(h)}$ denote the MLP embedding functions in layer $h$.
}

\rev{Finally, the output module maps the embedded graph obtained from the convolution module to scores of vessel pairs. Specifically,} the vessel nodes, whose embeddings now encode (after the convolutions) information about the whole graph, are input to an MLP representing a score for that pair. Finally, a softmax normalizes these scalars into a distribution over these pairs of vessels. 
%% \tr{RC: Is there any way to visualize this?}

\subsubsection{Loss Functions}
\label{sec_2_loss_function}
The simplest strategy is to treat the problem as a classification problem, where vessel pairs that had the highest expert score (could have been selected by the expert) are given the label 1, and the others are given the label 0. We then train with a cross-entropy loss to increase the probability of those pairs with label 1, and reduce the probability of those pairs with label 0.

Alternatives include giving label 1 to every vessel pair that has top 5 best expert scores \rev{(i.e., root bounds)}, instead of top 1. We have experimented with both two labeling strategies and the results of our performance are reported in Section \ref{sec_4_labelling}.

\subsubsection{Data Collection}
\label{sec:data}

As we have a limited number of instances available for training, it seems tricky to learn anything from a small set of instances at first glance. However, note that the neighborhood decisions have to be made multiple times during the solution process. Therefore, we can collect multiple samples in one run. Furthermore, we can apply the following procedure to generate multiple (potentially different) solution paths for each instance:\begin{enumerate}
\item At the beginning of each local search step, we flip a biased coin (with a random number generator) with the head being ``expert" (80$\%$ chance) and the tail being ``random" (20$\%$ chance).
\item If the outcome is the head (i.e., ``expert"), then we apply the expert \emph{Root} neighborhood selection rule, run the local search with the neighborhood selected by the expert, and collect the expert sample.
\item Otherwise, we randomly select the neighborhood and run the local search with the random neighborhood.
\item We then continue with the local search.
\end{enumerate}
With different random seeds, we run local searches on each instance multiple times following the hybrid neighborhood decision rules mentioned above. By doing so, we can collect even more data samples as (potentially) different sequences of samples can be collected within different runs of the same instance.

\subsubsection{The DAGGER Algorithm}
After collecting the initial training dataset following the approach described in Section \ref{sec:data}, the next step is to design the algorithm for training the GNN models on the collected dataset. The most straightforward way is to define the training problem as a one-time supervised learning problem, where the training data are collected offline before training and fixed forever. Then the GNNs will be trained by minimizing the supervised loss presented in Section \ref{sec_2_loss_function} on the training dataset.

A different line of research has been focusing on improving the out-of-sample performance of the GNN model. As our problem is essentially an online learning problem (decisions are made in a sequential manner), we \rev{apply} the DAGGER algorithm \citep{ross2011reduction} to train our GNN model. Specifically, the GNN model is trained iteratively using aggregated data that is collected from each iteration. In the first iteration, the data is generated from the trajectory given by the \emph{Root} policy ($\alpha\%$ chance following \emph{Root} and $(100-\alpha)\%$ chance selecting a random neighborhood for $\alpha\in\{60,70,80,90,100\}$). In the subsequent iterations, the data is generated by roughly following the policy associated with the GNN trained in the last iteration ($\alpha\%$ chance following GNN and $(100-\alpha)\%$ selecting a random neighborhood for $\alpha\in\{60,70,80,90,100\}$). As the algorithm requires several rounds of data collection and training, it is more time-consuming than the vanilla one-time approach.

\section{\rev{A Two-Phase Approach}}
\label{sec_algo}

The instances we consider can be categorized into two classes based on the constraints they include. The first class, which we refer to as combinatorial instances, consists of instances where only combinatorial constraints are present. 
\rev{These correspond to settings without tankfarm infrastructure, where each sales contract delivery is directly matched with a purchase contract lift. As a result, there are no inventory constraints involved.}
% \tr{RC: Maybe add some explanation of why we have combinatorial instances (some instances don't have tankfarms and all sales contracts are fulfilled by purchase contracts).} \rev{Nan, section 2.4, page 10: Combinatorial (Port Capacity and Ratability) Constraints}
The second class, called hybrid instances, includes both inventory and combinatorial constraints. For combinatorial instances (i.e., the first class), we apply a local search method, \rev{which is detailed in this section.} Our objective is to find a feasible solution to Problem \eqref{fea_MIP} such that $\sum_{i\in I}s_i=0$ in the objective function of Problem \eqref{fea_MIP}. 
For hybrid instances (the second class), we take a two-phase approach. In the first phase, we initially ignore the inventory constraints and focus solely on the combinatorial constraints. We first perform the local search to obtain the best solution under these constraints. Once this preliminary solution is found, in the second phase, we refine it by incorporating the inventory constraints and using off-the-shelf solvers to further optimize the solution.

\subsection{Local Search for Combinatorial Constraints}

For both combinatorial and hybrid instances, we rely on a local search algorithm to improve the solution \rev{by iteratively exploring its neighborhoods}. The details of the algorithm are provided in local search Algorithm~\ref{alg_ls}.

\begin{algorithm}[H]
\caption{Local Search for Combinatorial Constraints}
\label{alg_ls}
\begin{algorithmic}[1]
\State \textbf{Initialize} \rev{a feasible} solution $\vx^{\text{best}}$ \rev{of Problem \eqref{fea_MIP}} and its objective \rev{value} $z^{\text{best}}$ 
\While{Prespecified Termination Requirement}
\State \textbf{Select} two vessels $v, v' \in \V$ (a two‐vessel neighborhood) that have not been used since the last improvement based on the \rev{GNN for neighborhood search methods in Section~\ref{sec_neighborhood_gnn}}
\State \textbf{Fix} all binary \rev{vessel} variables \rev{(no slack variables are fixed)} not associated with $v$ or $v'$ (i.e., restrict the search to those indices in the chosen neighborhood)
\State \textbf{Solve} the resulting MIP subproblem to obtain a candidate solution $\vx^{\text{cand}}$ \rev{and its corresponding objective value $z(\vx^{\text{cand}})$}
\If{$z(\vx^{\text{cand}}) < z^{\text{best}}$} 
\State $\vx^{\text{best}} \gets \vx^{\text{cand}}$ and $z^{\text{best}} \gets z(\vx^{\text{cand}})$
\EndIf
\EndWhile
\State \textbf{Return} the final best solution $\vx^{\text{best}}$ and its objective $z^{\text{best}}$
\end{algorithmic}
\end{algorithm}
% \tr{RC: Is $z(\cdot)$ not defined?}

Local search is particularly effective for combinatorial instances, where the neighborhood structure plays a critical role in the effectiveness of local search. In our implementation, we focus on a two-vessel neighborhood, where each iteration modifies the decision variables associated with a selected pair of vessels while keeping all other variables fixed. \rev{We set the termination requirement in Step 2 of local search Algorithm~\ref{alg_ls} as a maximum number of $50$ iterations. At each iteration, we also set ``\texttt{NodeLimit}" to $1,000$, and ``\texttt{TimeLimit}" to $1,000$ seconds.} 
% \rev{Other} implementation details of \rev{local search} Algorithm~\ref{alg_ls} will be introduced in Section~\ref{sec_numerical}, 
This structured search of Algorithm~\ref{alg_ls} ensures that the problem remains computationally manageable while still exploring a diverse set of solutions. As discussed in Section~\ref{sec_ml_model}, the selection of vessel pairs (i.e., Step 3 in \rev{local search} Algorithm~\ref{alg_ls}) follows a systematic procedure guided by GNNs.

For hybrid instances, local search serves as the first phase of our solution approach. Given that these instances also involve inventory constraints, we initially disregard them, i.e., their violation is not penalized, and apply \rev{local search} Algorithm~\ref{alg_ls} to optimize the remaining optimization problem considering only combinatorial constraints. This allows us to efficiently generate a high-quality solution before proceeding to the refinement phase, where inventory constraints are incorporated for further optimization.

\subsection{Feasibility Refinement via Gurobi}

    After obtaining a high-quality solution from the local search phase, we refine it by incorporating inventory constraints and solving the full MIP model using Gurobi. 
    \rev{This refinement phase guarantees that the final solution to Problem~\eqref{fea_MIP} is feasible, with all slack variables equal to $0$ in both the combinatorial and inventory constraints.}
    % This refinement phase ensures that the final solution is feasible to Problem~\eqref{fea_MIP} with both combinatorial constraints and inventory constraints. 
    Naturally, since the local search phase (phase 1) only considers combinatorial constraints, the resulting solution may not fully satisfy inventory constraints. To address this, we initialize the MIP model that \rev{reinstalls} these inventory constraints with the best solution from local search (i.e., the output of \rev{local search} Algorithm~\ref{alg_ls}) as \rev{the initial solution} in Gurobi. The refinement phase is outlined in \rev{feasibility refinement} Algorithm~\ref{alg_feasibility_refinement}, which takes the local search output, integrates inventory constraints, and returns a refined solution after solving the full Problem \eqref{fea_MIP}. We remark that in instances with both combinatorial and inventory constraints, the slack variable $\bm{s}$ in Problem~\eqref{fea_MIP} is introduced for both types of constraints \rev{in the refinement phase} \rev{and the objective function of Problem~\eqref{fea_MIP} incorporates $\bm{s}$ for both the combinatorial and inventory constraints.} Accordingly, the set $I$ in Problem~\eqref{fea_MIP} includes both combinatorial and inventory constraints. 
\rev{Our preliminary experiments show that directly solving the full MIP problem without this two-phase approach incurs prohibitively long solution time, whereas our two‐phase approach reduces solution time dramatically.} 

\begin{algorithm}[H]
\caption{Feasibility Refinement via Gurobi}
\label{alg_feasibility_refinement}
\begin{algorithmic}[1]
\State \textbf{Input:} Solution $\vx^{\text{best}}$ from \rev{local search} Algorithm~\ref{alg_ls}
\State \rev{\textbf{Build}} MIP model \eqref{fea_MIP} with \rev{reinstalled} inventory constraints
\State \textbf{Set Initial Solution:} Use $\vx^{\text{best}}$ as \rev{the initial solution} for Gurobi
\State \textbf{Return} the final refined solution $\vx^{\text{best}}$
\end{algorithmic}
\end{algorithm}

\section{Computational Results}
\label{sec_numerical}

\rev{We conducted an extensive set of computational experiments to evaluate the performance of our proposed approaches.} All experiments were run on a computer with an Intel(R) Xeon(R) Gold 6258R processor running at 2.7 GHz, with up to four threads used. All optimization problems were solved using the optimization solver Gurobi with version 12.0.0. \rev{For each instance, we start our local search methods from the same solution of Problem \eqref{fea_MIP}, obtained by running Gurobi with ``\texttt{SolutionLimit}'' set to 1 to obtain an initial feasible solution and refining it using a local search on  the``single-vessel'' neighborhood for each vessel $v\in \V\cup\{0\}$.}

\subsection{Description of Instances}
We evaluate the performance of the proposed algorithm on a set of synthetic problem instances designed to closely resemble real-world LNG-IRO problems faced by a major vertically integrated energy company. In this \rev{section, we provide a} high-level overview of these instances to give the reader an understanding of their size and complexity (refer to Table~\ref{tab_problem_statistics} for a summary of the key statistics of the problem instances). 

Each instance is characterized primarily by the total LNG volume transported from supply locations to demand locations over a 365-day planning horizon. In this study, we vary this total volume between 8 and 20 million metric tonnes\rev{, consistent with portfolios managed by large LNG suppliers.}

Once the total volume for an instance is determined, we allocate it across two types of supply sources: 1) Third-party buying contracts, and 2) Tankfarms. Here, each instance is randomly assigned between 0 and 2 tankfarms, where inventory levels must be maintained within specified limits. On the demand side, we generate a set of sell contracts, ensuring that the total demand is slightly lower than the total supply. Any excess supply beyond contractual obligations can be sold through the spot market.

For each buy and sell contract, we impose constraints on the minimum and maximum number of deliveries or purchases required over different time intervals (e.g., 1, 3, 6, and 12 months).

Each contract and tankfarm is associated with a specific port that serves as the physical location for LNG loading and discharging. Ports have a predetermined number of berths, limiting the number of vessels that can be serviced simultaneously. We randomly assign real-world geographic locations to these ports \rev{to align with actual LNG supply and demand centers}. Sailing times between ports range from less than a day to over 65 days.

We generate vessels with capacities ranging from 155,000 to 175,000 cubic meters\rev{, which are representative of existing LNG vessels. This, combined with the portfolio volumes above, yields cargo counts consistent with real-world operations}. These instances include up to 25 vessels, with each of them randomly assigned one of two fuel modes: heavy fuel oil, which requires additional bunkering for fuel replenishment, and Boil-off LNG, which uses a portion of the transported LNG as fuel, eliminating the need for separate bunkering.

\rev{The generated instances are of realistic scale and capture the full range of operational constraints encountered in LNG MIRO, including supply-demand balancing, contractual obligations, port and berth limitations, sailing times, and vessel fuel considerations. To justify the use of synthetic instances and to make their scale explicit, we benchmark them against established references. These instances are comparable to, and in many cases larger than, those in seminal LNG inventory routing studies such as \cite{goel2012large}, and employ a 365-day planning horizon which is the maximum in the Group~2 class of problems in MIRPLib \citep{papageorgiou2014mirplib}. They are further complicated by the inclusion of multiple production ports alongside multiple consumption sites, increasing the combinatorial complexity. As shown in Table~\ref{tab_some_instances_comparisons}, the largest cases involve up to 25 vessels, 21 ports, and over 2.6~million variables, constituting a demanding test bed for the advanced heuristics proposed in this paper. The baseline performance of standard methods on these instances, also reported in Table~\ref{tab_some_instances_comparisons}, underscores their suitability for assessing algorithmic performance.}

\begin{table}[htbp]
\centering
\caption{Statistics of problem instances}
\renewcommand{\arraystretch}{1}
\setlength{\tabcolsep}{1pt}
\label{tab_problem_statistics}
\begin{tabular}{rrrr}
\hline
& Median & Min & Max \\ \hline
\# of vessels & 10 & 4 & 25 \\ 
\# of loading ports & 3 & 2 & 5 \\ 
\# of discharging ports & 9 & 5 & 16 \\ 
\# of continuous variables & 7,701 & 4,512 & 13,693 \\
\# of binary variables & 636,153 & 172,984 & 2,645,654 \\
\# of constraints & 143,827 & 44,777 & 406,955 \\ \hline
\end{tabular}
\end{table}

Table~\ref{tab_problem_statistics} presents the statistical summary of the problem instances considered in this paper, which provides an overview of the problem scale and complexity.
Table~\ref{tab_problem_statistics} includes the minimum, maximum, and median values for the number of vessels, loading ports, discharging ports, continuous variables, binary variables, and constraints across all instances.

\subsection{Phase 1: Instances with Combinatorial Constraints}

Recall that the instances are divided into two categories: (i) Those with only combinatorial constraints, and (ii) Those with both combinatorial and inventory constraints. For all instances, we first run \rev{local search}    Algorithm~\ref{alg_ls} in phase 1. For instances that include inventory constraints, the numerical results of solution refinement in phase 2 will be presented in Section \ref{sec_phase2}. Here, we evaluate the performance of the local search in Algorithm~\ref{alg_ls}. At each iteration, we set ``\texttt{NodeLimit}" to $1,000$, and ``\texttt{TimeLimit}" to $1,000$ seconds. Consequently, Step 5 of \rev{local search}  Algorithm~\ref{alg_ls} terminates if any of these limits is reached or if the optimal solution is found. Additionally, we configure Gurobi with ``\texttt{Heuristics}" set to $0.95$ and ``\texttt{Threads}" set to $4$.

\subsubsection{GNN Policies and Training}
In Step 3 of \rev{local search}  Algorithm~\ref{alg_ls}, we use different GNN policies to select the two-vessel neighborhood, with each policy named according to its training approach and implementation strategy.

\paragraph{Training Data Generation and Labeling Strategies}
\label{sec_4_labelling}

\begin{figure}
\includegraphics[width=\textwidth]{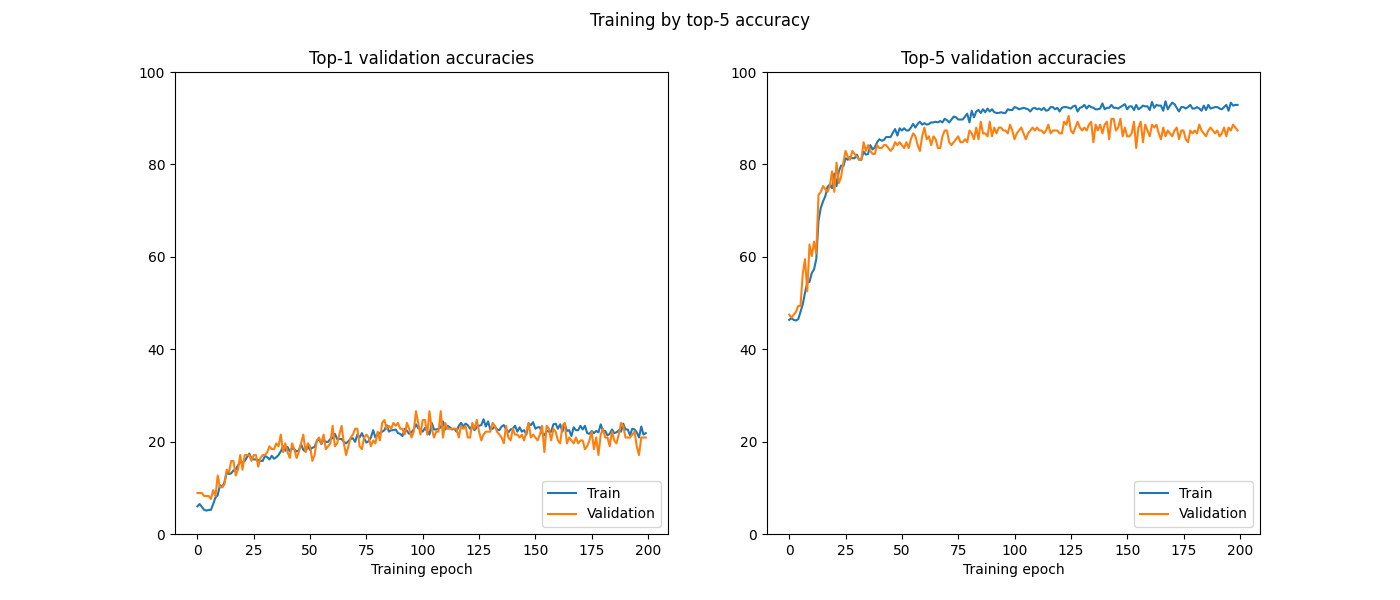}
\caption{Comparison of labeling strategies, when training with either top-1 (left) or top-5 (right) classification. The horizontal axes represent training epochs, and the vertical axes represent accuracy (in \%). The blue curves plot the training accuracies and the orange curves plot the validation accuracies of the implemented strategy accordingly.}
\label{fig:exp1-top1-top5}
\end{figure}

\rev{As presented in Section \ref{sec:data}, we have designed a data collection approach to augment training data by generating multiple solution paths from the same LNG-IRO instance. We conducted the data generation process on 72 available instances and collected a dataset of 5,502 training samples.} In order to label the training data, we considered two strategies as introduced in Section \ref{sec_2_loss_function}, giving label 1 to vessel pairs that have either top-1 (all vessel pairs with the best scores are targets)  or top-5 (all vessel pairs with scores in the top-5 best scores are targets) best expert scores. Specifically, we trained for 100 epochs, both using top-1 and top-5 classification. The results are reported in Figures \ref{fig:exp1-top1-top5}. As it can be seen, top-5 training gives uniformly better results, and with it, we can achieve roughly 20\% top-1 (i.e., standard) accuracy, and 80\% top-5 accuracy. However, we also observe a slight overfitting phenomenon: after a certain number of epochs the test performance plateaus or decreases, while the training accuracy keeps increasing. Note that we have implemented the $\ell_2$ regularization \citep{hinton2012practical} to mitigate overfitting. Based on our observation from the comparison, top-5 is finally selected for labeling the training data, and all the GNN models used in the subsequent experiments are trained using this strategy.

% \paragraph{How GNNs are implemented.} Note that given input, the GNN will output its estimate of how likely the neighborhood will be.  There are two particular ways that we implement our GNN models, namely, \textbf{GNN\_deter} and \textbf{GNN\_samp}. Specifically, \textbf{GNN\_deter} is a deterministic policy that selects the neighborhood with the highest score estimated by the particular GNN, while \textbf{GNN\_samp} is a sampling policy which samples a neighborhood from the distribution (over all possible neighborhoods) constructed using the GNN scores for all neighborhoods.

\paragraph{GNN Training.} We have trained and compared 4 GNN variants through the DAGGER training algorithm based on how the training data are split and the number of DAGGER training iterations. More specifically, the \textbf{GNN\_instance} is trained using the data split at the LNG instance level (80\% instances for training, and 20\% instances for evaluation), whereas \textbf{GNN\_sample} is trained with data split at local search neighborhood level (80\% neighborhood samples for training, and 20\% samples for evaluation). Under this scenario, the collected samples of vessel-pair neighborhoods from all the instances are mixed. We are interested in those two versions because we want to investigate how sensitive are the GNNs to the training data splitting at different granularity.  On the other hand, we also compare the GNN models trained by the DAGGER algorithm at different iteration levels. \textbf{GNN\_instance} and \textbf{GNN\_sample} are the GNN trained by only one iteration of DAGGER training, whereas \textbf{GNN\_instance2} and \textbf{GNN\_sample2} are the variants trained by two iterations of DAGGER training, respectively. That segmentation is to show how GNN behaves at different levels of training in terms of data intensity.

For notation convenience, in the following discussions, we denote \textbf{GNN\_instance}, \textbf{GNN\_instance2}, \textbf{GNN\_sample}, \textbf{GNN\_sample2} as $GNN\_i$, $GNN\_i2$, $GNN\_s$, $GNN\_s2$, respectively.

\subsubsection{Phase 1 Results}

\rev{Local search} Algorithm~\ref{alg_ls} iterates until no further objective improvement is attained through local search. We report the geometric mean across all instances on different methods. The geometric mean is used instead of the arithmetic mean as it better captures the central tendency of values that span multiple orders of magnitude, preventing dominance by extreme values \citep{mcalister1879xiii}. 

Table~\ref{tab_all_instances_phase_1_geo_average} presents the results, detailing the average number of iterations and the average of total computing time. \rev{Notice that the computing time in Table~\ref{tab_all_instances_phase_1_geo_average} is for computing a solution that is feasible only for combinatorial constraints.} Since the final output of \rev{local search}  Algorithm~\ref{alg_ls} may be positive, i.e., $\sum_{i\in I}s_i>0$, this indicates that a feasible solution was not found. In such cases, we classify the instance as non-convergent. The table reports the number of non-convergent instances for each method.
Table~\ref{tab_all_instances_phase_1_geo_average} provides separate statistics for (i) All instances combined, (ii) \rev{Combinatorial} instances with combinatorial constraints only, and (iii) \rev{Hybrid} instances with combinatorial constraints and inventory constraints. 

From the results, we observe that methods based on instance-level GNN models ($GNN\_i$ and $GNN\_i2$) generally require fewer iterations compared to their sample-based counterparts ($GNN\_s$ and $GNN\_s2$).

Among all the methods, $GNN\_i$ emerges as the most efficient approach, achieving the lowest computing time. Specifically, $GNN\_i$ attains an average runtime of $\rev{2,456.19}$ seconds across all instances, which is lower than the $\rev{2,804.58}$ seconds required by $GNN\_s$ and the $\rev{2,496.65}$ seconds by $GNN\_s2$. Additionally, $GNN\_i$ requires fewer iterations on average $(10.26)$ compared to $GNN\_s$ $(12.24)$ and $GNN\_s2$ $(12.75)$, further demonstrating its effectiveness in rapidly converging to high-quality solutions. More importantly, for all instances containing only combinatorial constraints, $GNN\_i$ successfully provides a feasible solution \rev{with all combinatorial constraints satisfied}, while $GNN\_s$ does not. \rev{For instances involving only combinatorial constraints, we apply the model-based selection strategy of \cite{goel2012large} for comparison purposes. The main idea of the model-based selection strategy of \cite{goel2012large} is to roughly identify the two-vessel neighborhood subproblem (by LP relaxation) that will lead to the biggest improvement in the objective function with respect to the current solution. However, because our problem is more complex than the model tested in \cite{goel2012large}, its performance is not competitive. Out of the 37 test instances, the model-based selection strategy produces a feasible solution to Problem~\eqref{fea_MIP} in only 15 cases; for the other 22 instances, it fails to return any feasible solution. For most of the failed cases, the reason of failure is that the neighborhood selection MIP model cannot be solved within $1,000$ seconds.} \rev{We also implement the random heuristic of \cite{goel2012large} and the result is presented in the last row of Table~\ref{tab_all_instances_phase_1_geo_average}, where the main idea is to select the two-vessel neighborhood subproblem randomly. It is computationally cheap per iteration, though it requires more iterations overall. However, this heuristic does not always produce a feasible solution for Problem~\eqref{fea_MIP} across all instances. In phase 2, we will also evaluate and compare its solution quality with our $GNN\_i$ method. Note that in phase 1 we are only looking for solutions that satisfy all combinatorial constraints.}

\begin{table}[htbp]
\centering
\caption{Geometric mean across all instances}
\renewcommand{\arraystretch}{1.5}
\setlength{\tabcolsep}{6pt}
\label{tab_all_instances_phase_1_geo_average}
\rev{
\resizebox{\textwidth}{!}{
\begin{tabular}{cccccccccc}
\hline
\multirow{2}{*}{Methods} & \multicolumn{3}{c}{Average of All Instances}                  & \multicolumn{3}{c}{\makecell{Average of Instances With \\Only Combinatorial Constraints}}                & \multicolumn{3}{c}{\makecell{Average of Instances With \\ Inventory Constraints}}                                         \\ \cline{2-10} 
& \multicolumn{1}{c}{Iterations} & \multicolumn{1}{c}{Time (s)} & \multicolumn{1}{c}{\makecell{\# of Non-convergent  \\(\# Total) Instances}} & \multicolumn{1}{c}{Iterations} & \multicolumn{1}{c}{Time (s)} & \multicolumn{1}{c}{\makecell{\# of Non-convergent \\ (\# Total) Instances}} & \multicolumn{1}{c}{Iterations} & \multicolumn{1}{c}{Time (s)}  & \multicolumn{1}{c}{\makecell{\# of Non-convergent\\ (\# Total) Instances}} \\ \hline
$GNN\_i$            & \multicolumn{1}{r}{10.26}      & \multicolumn{1}{r}{2,456.19} & 2 (57)                                  & \multicolumn{1}{r}{6.31}       & \multicolumn{1}{r}{1,544.11}  & 0 (37)                                  & \multicolumn{1}{r}{22.45}      & \multicolumn{1}{r}{5,182.64}  & 2 (20)                                  \\ \hline
$GNN\_i2$           & \multicolumn{1}{r}{10.86}      & \multicolumn{1}{r}{2,544.05} & 4 (57)                                  & \multicolumn{1}{r}{6.30}       & \multicolumn{1}{r}{1,518.38}  & 0 (37)                                  & \multicolumn{1}{r}{26.09}      & \multicolumn{1}{r}{5,835.88} & 4 (20)                                  \\ \hline
$GNN\_s$              & \multicolumn{1}{r}{12.24}      & \multicolumn{1}{r}{2,804.58} & 7 (57)                                  & \multicolumn{1}{r}{7.56}       & \multicolumn{1}{r}{1,845.39} & 1 (37)                                  & \multicolumn{1}{r}{26.57}      & \multicolumn{1}{r}{5,499.14} & 6 (20)                                  \\ \hline
$GNN\_s2$           & \multicolumn{1}{r}{12.75}      & \multicolumn{1}{r}{2,496.65} & 6 (57)                                  & \multicolumn{1}{r}{8.30}       & \multicolumn{1}{r}{1,623.31} & 3 (37)                                  & \multicolumn{1}{r}{25.41}      & \multicolumn{1}{r}{4,990.12}  & 3 (20)                                  \\ \hline
$\mathrm{Random}$           & \multicolumn{1}{r}{13.88}      & \multicolumn{1}{r}{3,633.77} & 9 (57)                                  & \multicolumn{1}{r}{9.14}       & \multicolumn{1}{r}{2,149.95} & 4 (37)                                  & \multicolumn{1}{r}{31.33}      & \multicolumn{1}{r}{ 10,097.64}  & 5 (20)                                  \\ \hline
\end{tabular}
}
}
\end{table}

We also compare our proposed method with solving the instances using Gurobi directly. Table~\ref{tab_gurobi_vs_gnn_instance_without_inventory} presents the geometric mean of computing time for instances with only combinatorial constraints and the results demonstrate a significant improvement in computational efficiency. \rev{$GNN\_i$} achieves an average runtime of \rev{$2,456.19$} seconds, compared to $14,231.76$ seconds for Gurobi, resulting in a \rev{82.74\%} reduction in computing time.  \rev{For all instances containing only combinatorial constraints, our method $GNN\_i$ successfully produces a feasible solution to Problem~\eqref{fea_MIP}. In these cases, the objective values obtained by the $GNN\_i$ method match those produced by Gurobi within the longer runtime, indicating that our approach not only finds feasible solutions but also maintains competitive solution quality.}
% This improvement suggests that our approach effectively accelerates the solution process while maintaining solution quality.
Notice that, in Table~\ref{tab_all_instances_phase_1_geo_average}, $GNN\_i2$ outperforms $GNN\_i$ for instances with only combinatorial constraints in terms of both average running time and average number of iterations. However, when considering all instances, $GNN\_i$ demonstrates superior overall performance.

\begin{table}[htbp]
\centering
\caption{\rev{Computing time comparison} between Gurobi and $GNN\_i$ on instances with only combinatorial constraints}
\renewcommand{\arraystretch}{1}
\setlength{\tabcolsep}{2pt} % Adjust column spacing
\label{tab_gurobi_vs_gnn_instance_without_inventory}
\begin{tabular}{cr}
\hline
\textbf{Methods}      & \textbf{Time (s)} \\ \hline
Gurobi                & 14,231.76                 \\ \hline
\rev{$GNN\_i$}       & \rev{2,456.19}                  \\ \hline
 Improvement           & \rev{82.74}\%                  \\ \hline
\end{tabular}
\end{table}

\rev{Finally, Table~\ref{tab_some_instances_comparisons} reports results for a selection of representative instances chosen to illustrate the range of problem sizes considered in this paper. The set includes three relatively small (``easy”), three medium, and four large (``hard”) instances, depending on the total number of variables. These sizes are substantially larger than those typically studied in the literature (see, e.g., \cite{goel2012large}). The table presents both the instance-specific statistics (number of vessels, loading ports, and discharging ports) and the runtime of different solution methods. The results highlight that our proposed method maintains competitive performance across all sizes and scales effectively for the largest instances. In contrast, Gurobi, the model-based selection strategy of \cite{goel2012large}, and the random heuristic of \cite{goel2012large} exhibit significant performance degradation on larger instances.}

\begin{table}[htbp]
\centering
\caption{\rev{Computing comparison on the representative instances}}
\renewcommand{\arraystretch}{1}
\setlength{\tabcolsep}{1pt}
\label{tab_some_instances_comparisons}
\rev{
\resizebox{\textwidth}{!}{
\begin{tabular}{ccrcccrrrrrrr}
\hline
\multicolumn{2}{c}{\multirow{2}{*}{Instances}}   & \multicolumn{4}{c}{\textbf{Model Parameters}}                   & \multicolumn{7}{c}{\textbf{Comparisons of Methods by Time (s)}}                \\ \cline{3-13} 
\multicolumn{2}{c}{}                  & \multicolumn{1}{c}{\makecell{\# of Total \\  Variables}} & \multicolumn{1}{c}{\makecell{\# of \\ Vessels}} & \multicolumn{1}{c}{\makecell{\# of Loading\\  Ports}} & \multicolumn{1}{c}{\makecell{\# of Discharging\\  Ports}} & \multicolumn{1}{c}{$GNN\_i$} & \multicolumn{1}{c}{$GNN\_i2$} & \multicolumn{1}{c}{$GNN\_s$} & \multicolumn{1}{c}{$GNN\_s2$} & \multicolumn{1}{c}{Gurobi}    & \multicolumn{1}{c}{\makecell{Model-based\\Selection}} & \multicolumn{1}{c}{\makecell{Random\\Heuristic}} \\ \hline
\multicolumn{1}{c}{\multirow{3}{*}{easy}}   & \phantom{0}1  & \multicolumn{1}{r}{177,496}           & {\phantom{0}5}             & {\phantom{0}2}                   & \phantom{0}5                                            & \multicolumn{1}{r}{81.88}    & \multicolumn{1}{r}{79.32}     & \multicolumn{1}{r}{81.25}    & \multicolumn{1}{r}{82.18}     & \multicolumn{1}{r}{2,252.84}   & \multicolumn{1}{r}{56.19}       & 57.81                       \\ \cline{2-13} 
\multicolumn{1}{c}{}                        & \phantom{0}2  & \multicolumn{1}{r}{235,678}           & \multicolumn{1}{c}{\phantom{0}7}             & \multicolumn{1}{c}{\phantom{0}2}                   & \phantom{0}5                                            & \multicolumn{1}{r}{97.53}    & \multicolumn{1}{r}{107.43}    & \multicolumn{1}{r}{93.20}    & \multicolumn{1}{r}{110.71}    & \multicolumn{1}{r}{2,790.82}   & \multicolumn{1}{r}{94.38}       & 80.29                       \\ \cline{2-13} 
\multicolumn{1}{c}{}                        & \phantom{0}3  & \multicolumn{1}{r}{238,969}           & \multicolumn{1}{c}{\phantom{0}7}             & \multicolumn{1}{c}{\phantom{0}2}                   & \phantom{0}5                                            & \multicolumn{1}{r}{108.76}   & \multicolumn{1}{r}{118.87}    & \multicolumn{1}{r}{117.26}   & \multicolumn{1}{r}{119.08}    & \multicolumn{1}{r}{4,352.92}   & \multicolumn{1}{r}{92.55}       & 82.38                       \\ \hline
\multicolumn{1}{c}{\multirow{3}{*}{medium}} & \phantom{0}4  & \multicolumn{1}{r}{589,710}           & \multicolumn{1}{c}{12}            & \multicolumn{1}{c}{\phantom{0}3}                   & \phantom{0}5                                            & \multicolumn{1}{r}{1,620.43}  & \multicolumn{1}{r}{1,120.80}   & \multicolumn{1}{r}{1,348.13}  & \multicolumn{1}{r}{1,530.98}   & \multicolumn{1}{r}{10,276.68}  & \multicolumn{1}{r}{5,836.18}     & 2,573.33                     \\ \cline{2-13} 
\multicolumn{1}{c}{}                        & \phantom{0}5  & \multicolumn{1}{r}{599,218}           & \multicolumn{1}{c}{\phantom{0}5}             & \multicolumn{1}{c}{\phantom{0}3}                   & \phantom{0}9                                            & \multicolumn{1}{r}{4,665.35}  & \multicolumn{1}{r}{5,559.00}   & \multicolumn{1}{r}{12,590.37} & \multicolumn{1}{r}{4,248.70}   & \multicolumn{1}{r}{11,934.00}  & \multicolumn{1}{r}{$\mathrm{Time Limit}$}       & 3,301.45                     \\ \cline{2-13} 
\multicolumn{1}{c}{}                        & \phantom{0}6  & \multicolumn{1}{r}{602,142}           & \multicolumn{1}{c}{\phantom{0}8}             & \multicolumn{1}{c}{\phantom{0}2}                   & 13                                           & \multicolumn{1}{r}{1,823.20}  & \multicolumn{1}{r}{1,928.75}   & \multicolumn{1}{r}{2,220.38}  & \multicolumn{1}{r}{1,832.21}   & \multicolumn{1}{r}{8,809.17}   & \multicolumn{1}{r}{99,898.88}    & 3,133.64                     \\ \hline
\multicolumn{1}{c}{\multirow{4}{*}{hard}}   & \phantom{0}7  & \multicolumn{1}{r}{1,698,804}          & \multicolumn{1}{c}{10}            & \multicolumn{1}{c}{\phantom{0}4}                   & 11                                           & \multicolumn{1}{r}{7,578.76}  & \multicolumn{1}{r}{7,904.95}   & \multicolumn{1}{r}{11,090.28} & \multicolumn{1}{r}{8405.70}   & \multicolumn{1}{r}{138,399.90} & \multicolumn{1}{r}{106,922.21}   & 27,554.75                    \\ \cline{2-13} 
\multicolumn{1}{c}{}                        & \phantom{0}8  & \multicolumn{1}{r}{1,722,572}          & \multicolumn{1}{c}{11}            & \multicolumn{1}{c}{\phantom{0}3}                   & 13                                           & \multicolumn{1}{r}{7,661.17}  & \multicolumn{1}{r}{9,919.45}   & \multicolumn{1}{r}{10,341.45} & \multicolumn{1}{r}{8,181.08}   & \multicolumn{1}{r}{59,836.42}  & \multicolumn{1}{r}{106,548.92}   & 14,442.40                    \\ \cline{2-13} 
\multicolumn{1}{c}{}                        & \phantom{0}9  & \multicolumn{1}{r}{2,273,522}          & \multicolumn{1}{c}{10}            & \multicolumn{1}{c}{\phantom{0}4}                   & 16                                           & \multicolumn{1}{r}{11,406.09} & \multicolumn{1}{r}{10,469.20}  & \multicolumn{1}{r}{11,363.12} & \multicolumn{1}{r}{13,286.45}  & \multicolumn{1}{r}{107,056.15} & \multicolumn{1}{r}{111,071.74}   & 37,846.52                    \\ \cline{2-13} 
\multicolumn{1}{c}{}                        & 10 & \multicolumn{1}{r}{2,659,347}          & \multicolumn{1}{c}{10}            & \multicolumn{1}{c}{\phantom{0}4}                   & 16                                           & \multicolumn{1}{r}{16,740.47} & \multicolumn{1}{r}{16,160.01}  & \multicolumn{1}{r}{14,546.42} & \multicolumn{1}{r}{18,611.62}  & \multicolumn{1}{r}{271,328.80} & \multicolumn{1}{r}{113,780.67}   & 38,486.50                    \\ \hline
\end{tabular}
}
}
\end{table}

\subsection{Phase 2: Results \rev{on Hybrid Instances}}
\label{sec_phase2}

For hybrid instances with both combinatorial constraints and inventory constraints, as described in \rev{feasibility refinement} Algorithm~\ref{alg_feasibility_refinement}, we use the solution obtained in phase 1 as \rev{the Gurobi initial solution} to enhance computational efficiency. Notice that, as shown in Table~\ref{tab_all_instances_phase_1_geo_average}, the solution obtained from \rev{local search}  Algorithm~\ref{alg_ls} may not always be feasible \rev{(i.e., satisfy all combinatorial constraints)}. Specifically, 2 out of 57 instances in Table~\ref{tab_all_instances_phase_1_geo_average} exhibit this issue with the $GNN\_i$ method. Despite this, we use the precomputed solution from Table~\ref{tab_all_instances_phase_1_geo_average} to initialize the hybrid instances. By doing so, we aim to accelerate the search process and guide Gurobi more effectively toward a feasible solution \rev{that satisfies all combinatorial and inventory constraints}. This warm-start strategy is particularly beneficial in large-scale instances where finding \rev{a good initial} solution can be computationally expensive. To assess the effectiveness of our approach, we compare it against a baseline where the instances with inventory constraints \rev{(modeled using Problem \eqref{fea_MIP})} are solved directly using Gurobi without \rev{an initial solution}. This comparison allows us to quantify the impact of our method on solution time. The results are presented in Table~\ref{tab_gurobi_vs_gnn_instance_with_inventory}, where we report the geometric average runtime of our method compared to solving the instances directly using Gurobi. The running time of our method includes both the total time spent on the local search in Algorithm~\ref{alg_ls} with the $GNN\_i$ method and the refinement step using Gurobi in \rev{feasibility refinement} Algorithm~\ref{alg_feasibility_refinement}. Our approach demonstrates a solid computational advantage, with our method helping Gurobi converge to \rev{optimality} faster.  \rev{We also benchmark our approach against the random heuristic from \cite{goel2012large}. While that heuristic is fast in phase 1, the lower-quality solutions it produces can worsen performance in phase 2, increasing its total runtime.}

\begin{table}[htbp]
\centering
\caption{\rev{Computing time comparisons among Gurobi, $\mathrm{Random}$, and $GNN\_i$ for hybrid instances}}
\renewcommand{\arraystretch}{1}
\setlength{\tabcolsep}{2pt} % Adjust column spacing
\label{tab_gurobi_vs_gnn_instance_with_inventory}
\rev{
\begin{tabular}{cr}
\hline
\textbf{Methods}      & \textbf{Time (s)}  
 % \textbf{Improvement (\%)} 
\\ \hline 
Gurobi                & 122,539.75       \\ \hline
$\mathrm{Random}$            & 51,120.30          \\ \hline
\rev{$GNN\_i$} & 14,589.33  \\
\hline
\end{tabular}
}
\end{table}

\section{Conclusion}
\label{sec_conclusion}

In this work, we proposed a machine learning-based local search approach for finding feasible solutions of large-scale maritime inventory routing optimization problems. Given the combinatorial complexity of MIRO, particularly in \rev{hybrid instances} with inventory constraints, we integrated a GNN-based neighborhood selection method to enhance local search efficiency. Our approach enables a structured exploration of vessel neighborhoods, improving solution quality while maintaining computational efficiency. Through extensive computational experiments on synthetic instances, we demonstrated that our method outperforms direct MIP approaches \rev{as well as existing local search approaches, achieving substantial improvements in both solution time and quality.}
% in solution time. 
Beyond its immediate application in MIRO, our methodology offers broader implications for learning-based optimization in large-scale industrial logistics, \rev{highlighting how the integration of optimization models with machine learning can address scalability challenges in complex systems.}
% problems. 
Future research directions include extending the GNN framework to handle additional problem structures, incorporating uncertainty in demand and transit times, and exploring different learning techniques for further improving neighborhood selection strategies. \rev{Overall, this work highlighted the potential of learning-based optimization as a powerful paradigm for tackling complex, data-driven challenges in maritime inventory routing.}

\bibliography{ref.bib}

\end{document}